\documentclass[twoside,reqno]{lt6-proc}
\usepackage{epsfig,cite}
\usepackage{amssymb,amsmath}
\usepackage{times}
\def\nn{\nonumber}
\def\beq{\begin{equation}}
\def\eeq{\end{equation}}
\def\bea{\begin{eqnarray}}
\def\eea{\end{eqnarray}}
\def\SL{SL_q(2)}
\def\U{{\cal U}}
\def\A{{\cal A}}
\def\e#1#2{e^{#1}_{#2}}
\def\CGC#1#2{C^{#1}_{#2}\,}
\def\D#1#2{D^{\ #1}_{#2}}
\def\T#1#2{T^{\ (#1)}_{#2}}
\def\h#1{\hat{#1}}
\def\pair#1#2{\left\langle #1,#2 \right\rangle}
\begin{document}
\sloppy \raggedbottom
\setcounter{page}{1}

\newpage
\setcounter{figure}{0}
\setcounter{equation}{0}
\setcounter{footnote}{0}
\setcounter{table}{0}
\setcounter{section}{0}



\title{Noncommutative Superspaces Covariant Under \boldmath{$OSp_q(1|2)$} Algebra}

\runningheads{N. Aizawa and R. Chakrabarti}{Noncommutative Covariant Spaces of $OSp_q(1|2)$}

\begin{start}


\author{N. Aizawa}{1},
\coauthor{R. Chakarabarti}{2},

\address{Department of Mathematics and Information Sciences, Graduate School of Science, 
Osaka Prefecture University, Daisen Campus, Sakai, Osaka 590-0035, Japan}{1}
\address{Department of Theoretical Physics, University of Madras, Guindy Campus, Chennai 
600 025, India}{2}


\begin{Abstract}
 Using the corepresentation of the quantum group $ SL_q(2)$ 
a general method for constructing noncommutative spaces 
covariant under its coaction is developed. 
The method allows us to treat the quantum plane and Podle\'s' 
quantum spheres in a unified way and to construct higher 
dimensional noncommutative spaces systematically. 
Furthermore, we extend the method to the quantum supergroup 
$ OSp_q(1|2).$ 
In particular, a one-parameter family of covariant algebras, 
which may be interpreted as noncommutative superspheres, is 
constructed. 
\end{Abstract}
\end{start}


\section{Introduction}

  Quantum groups provide a very powerful tool for 
investigations of noncommutative geometry, since 
they may be regarded as a noncommutative extension 
of linear Lie groups. Pioneering works by Manin\cite{Ma}, 
Woronowicz\cite{Wor}, Wess and Zumino\cite{WZ} 
are followed by hundreds of publications 
(see for example \cite{KS,BL} and references therein). 
One way of using quantum groups for noncommutative geometry 
is, by regarding them as an 
example of noncommutative manifolds, to develop a 
harmonic analysis on quantum groups. 
Another way, which may be more familiar to large class of 
physicists, is to regard a quantum group as a transformation 
matrix of vectors. Because of the noncommutative nature of 
quantum groups, the vectors transformed by a quantum group are 
{\it a priori} noncommutative. Namely, the components of the 
vectors have nontrivial commutation relations. 
In order to fit such vectors in theories of physics, 
we require covariance, that is, the commutation relations 
are preserved by quantum group transformations. 
Noncommutative vectors obeying a covariant algebra may be  
given geometrical interpretation. 

  Let us consider $ SU_q(2) $ as an example. 
If we consider a covariant algebra transformed by 
the fundamental representation, it may be interpreted as 
a noncommutative analogue of two dimensional flat space\cite{Ma}. 
While if we take a covariant algebra for 
adjoint representation, it may be regarded as a noncommutative 
extension of 3-sphere\cite{Po87}. 
Higher dimensional representations 
will give higher dimensional noncommutative spaces. 
However, no such work has been done because mainly of 
computational difficulty. 
Furthermore, there are only a few works on noncommutative 
analogues of superspaces in the context of quantum groups 
despite the fact that supersymmetry is one of the most 
important notions of theoretical physics. 

  In the present work, using the corepresentations of 
quantum groups, we develop a general method for 
constructing noncommutative spaces for the simplest 
and the most important quantum (super) groups 
$ SL_q(2) $ and  $ OSp_q(1|2).$ In the first part of 
this paper (\S \ref{SLcorep} and \S \ref{SL2cov}), 
the case $ SL_q(2)$ is considered. It will be 
seen that, by our method, the quantum plane and 
the quantum spheres are treated on the same footing 
and that the higher dimensional noncommutative spaces 
may be constructed systematically. In the second 
part (\S \ref{CovSpeOSP}), we extend the 
method to $ OSp_q(1|2).$ As an application of our 
method, noncommutative superspace and a one-parameter 
family of noncommutative superspheres are explicitly 
constructed. Finally \S \ref{ConcludingR} is devoted to 
concluding remarks.


\section{\boldmath{$\SL$} and its corepresentations}
\label{SLcorep}

This section is a brief review of the definitions and 
representations for the quantum group $ SL_q(2) $ and 
the quantum algebra $ U_q[sl(2)] $ that is dual to $ SL_q(2).$ 
There are several good textbooks on this topics. 
Readers may refer, for example, to \cite{KS,BL} and 
references therein. 

  The quantum group $\SL$ is generated by four 
elements $a,b,c $ and $d$ subject to the relations
\bea
  & & a b = q b a, \qquad a c = q c a, \qquad b d = q d b, 
  \nn \\
  & & c d = q d c, \qquad b c = c b, \qquad \ ad - da = (q-q^{-1}) bc,
  \label{SLq2Def} \\
  & & ad - qbc = da - q^{-1}bc = 1. \nn
\eea
As is well-known, the coproduct ($ \Delta$), 
the counit ($\epsilon$) and the antipode ($S$) defined as follows 
make $\SL$ a Hopf algebra:
\bea
 & & 
  \Delta
  \begin{pmatrix}
     a & b \\ c & d
  \end{pmatrix}
  =
  \begin{pmatrix}
     a & b \\ c & d
  \end{pmatrix}
  \stackrel{\cdot}{\otimes}
  \begin{pmatrix}
     a & b \\ c & d
  \end{pmatrix},
  \qquad
  \epsilon
  \begin{pmatrix}
     a & b \\ c & d
  \end{pmatrix}
  = 
  \begin{pmatrix}
    1 & 0 \\ 0 & 1
  \end{pmatrix},
 \nn \\
 & & 
  S
  \begin{pmatrix}
     a & b \\ c & d
  \end{pmatrix}
  =
  \begin{pmatrix}
    d & -q^{-1}b \\ -qc & a
  \end{pmatrix}.
  \label{HopfSL2}
\eea
With the Hopf algebra mappings, we define a corepresentation 
of a quantum group. 
A vector space $V$ is called a right $\SL$-comodule 
if there exists a linear mapping 
$ \varphi_R : V \rightarrow V \otimes \SL $ satisfying
\beq
  (\varphi_R \otimes {\rm id})\circ \varphi_R = 
  ({\rm id}\otimes \Delta) \circ \varphi_R,
  \qquad
  ({\rm id}\otimes \epsilon) \circ \varphi_R = {\rm id}.
  \label{Rcomodule}
\eeq
Similarly, the left $\SL$-comodule is defined as a vector 
space $V$ equipped with a linear mapping 
$ \varphi_L : V \rightarrow \SL \otimes V $ such that
\beq
  ({\rm id}\otimes \varphi_L) \circ \varphi_L = 
  (\Delta \otimes {\rm id}) \circ \varphi_L,
  \qquad
  (\epsilon \otimes {\rm id}) \circ \varphi_L = {\rm id}.
  \label{Lcomodule}
\eeq
The mapping $ \varphi_R $ ($ \varphi_L $) is called 
a corepresentation, or, equivalently, a right (left) coaction of $ \SL $ 
on $V$. 
It is known that each irreducible corepresentation of $\SL$ is, 
as classical $SL(2),$ 
specified by the highest weight $j$ which takes any nonnegative integral 
or half-integral values. 
Let $V^{(j)}$ be a right $\SL$-comodule with the highest 
weight $j$ and 
$ \{ e^j_m, \ m = j, j-1, \cdots , -j \} $ be its basis:
\beq
   \varphi_R(e^j_m) = \sum_{m'} e^j_{m'} \otimes T^{(j)}_{m'm},
   \quad
   T^{(j)}_{m'm} \in \SL
   \label{Raction}
\eeq
The corepresentations of $\SL $ have been obtained explicitly\cite{Mas,VS,Koo}. 
We here give $j=1/2$ and $ j=1$ corepresentation matrices as an example:
\beq
 T^{(1/2)} = 
 \begin{pmatrix}
   a & b \\ c & d
 \end{pmatrix},
 \label{thalfmatrix}
\eeq
\beq
  T^{(1)} = \left(
    \begin{array}{ccc}
      a^2 & (1+q^{-2})^{1/2} ab & b^2 \\
     (1+q^{-2})^{1/2}ac & 1 + (q+q^{-1})bc & (1+q^{-2})^{1/2} bd \\
     c^2 & (1+q^{-2})^{1/2}cd & d^2
    \end{array}
  \right).
  \label{t1matrix}
\eeq

 A comodule of a quantum group is, in general, a module {\it i.e.} 
a representation space of the dual quantum algebra. 
We define the action of $ U_q[sl(2)]$ on $V^{(j)}$ by
\beq
  X e^j_m = (({\rm id} \otimes X) \circ \varphi_R)(e^j_m)
  = \sum_{m'} e^j_{m'} \left\langle X, T^{(j)}_{m'm} \right\rangle,
  \quad X \in U_q[sl(2)]
  \label{actionofUonV}
\eeq
where $ \langle \ ,\ \rangle : U_q[sl(2)] \otimes \SL \rightarrow {\mathbb C} $ 
is the duality pairing of two Hopf algebras.  
Then it may be verified that 
the matrix $ \left\langle X, T^{(j)}_{m'm} \right\rangle $ gives 
an irreducible representation of $ U_q[sl(2)]$ with the highest weight $j$. 
The product space $ V^{(j_1)} \otimes V^{(j_2)} $ is, in general, reducible 
and is decomposed into irreducible spaces as
\beq
  j_1 \otimes j_2 = j_1 + j_2 \oplus j_1 + j_2 -1 \oplus \cdots 
  \oplus |j_1-j_2|.
  \label{Irrep_decomp} 
\eeq
The decomposition is carried out by the Clebsch-Gordan coefficients (CGC)
\beq
  e^J_M(j_1,j_2) = \sum_{m_1,m_2} C^{\;j_1j_2J}_{m_1m_2M} e^{j_1}_{m_1} 
  \otimes e^{j_2}_{m_2}.
  \label{j1j2toJ}
\eeq
The CGC satisfy the following orthogonality relations
\bea
  & & \sum_{j,m} C^{\;j_1j_2j}_{m_1m_2m} C^{\;j_1j_2j}_{m'_1m'_2m'} 
      = \delta_{m_1m'_1} \delta_{m_2m'_2},
  \label{CGCortho1} \\
  & & \sum_{m_1,m_2} C^{\;j_1j_2j}_{m_1m_2m} C^{\;j_1j_2j'}_{m_1m_2m'} 
      = \delta_{jj'} \delta_{mm'}.
  \label{CGCortho2}
\eea
Two corepresentations are also coupled by CGC. 
The coupling is given by the formula that is called 
the Wigner's product law\cite{GKK}:
\beq
  \delta_{jj'} T^{(j)}_{mm'} = \sum_{m_1, m_2 \atop m'_1,m'_2}
  C^{\;j_1j_2j}_{m_1m_2m} C^{\;j_1j_2j'}_{m'_1m'_2m'} 
  T^{(j_1)}_{m_1m'_1} T^{(j_2)}_{m_2m'_2}.
  \label{WigPro}
\eeq

%
\section{Covariant algebras of \boldmath{$\SL$}}
\label{SL2cov}

\subsection{General prescription}

In this section, we give a general prescription to 
construct $\SL$-covariant algebras. 
By covariant algebras, we mean algebras whose defining 
relations are preserved under the right coaction of $\SL$ 
defined by (\ref{Raction}). 
Probably, the simplest way to find such an algebra 
is to introduce an algebraic structure on the comodule $ V^{(j)}.$ 
Let $ \mu $ be a product in $ V^{(j)},\ i.e.,$ 
$ \mu(f \otimes g ) = fg, \ f, g \in V^{(j)}.$ 
We specifically consider the following composite object
\beq
  \mu(e^J_M(j,j)) = \sum_{m_1,m_2} C^{\; jjJ}_{m_1m_2M} 
  e^j_{m_1} \,e^j_{m_2}.
  \label{muCouple}
\eeq
The right coaction on (\ref{muCouple}) is shown to be
\beq
  \varphi_R \circ \mu(e^J_M(j,j)) = \sum_{M'} \mu(e^{J}_{M'}(j,j)) \otimes T^{(J)}_{M'M}.
  \label{covarmu}
\eeq
The proof may be done in a straightforward way by inverting the relation 
(\ref{muCouple})
\beq
   e^j_{m_1} e^j_{m_2} = \sum_{JM} C^{\;jjJ}_{m_1m_2M} \mu(e^J_M(j,j)),
   \label{inverseC}
\eeq
and subsequently using the product law (\ref{WigPro})
\bea
  & &  \varphi_R \circ \mu(e^J_M(j,j)) 
  \nn \\
  & &\quad = 
  \sum_{m_1m_2} C^{jjJ}_{m_1m_2M} \varphi_R(e^j_{m_1}) \varphi_R(e^j_{m_2})
  \nn \\
  & & \quad = 
   \sum_{m_1m_2 \atop m'_1m'_2} C^{jjJ}_{m_1m_2M} e^j_{m'_1}e^j_{m'_2}
   \otimes T^{(j)}_{m'_1m_1} T^{(j)}_{m'_2m_2}
   \nn \\
  & & \quad \stackrel{(\ref{inverseC})}{=}
    \sum_{m_1m_2 \atop m'_1m'_2} \sum_{J'M'} C^{jjJ}_{m_1m_2M} C^{jjJ'}_{m'_1m'_2M'}
      \mu(e^{J'}_{M'}(j,j)) 
  \otimes T^{(j)}_{m'_1m_1} T^{(j)}_{m'_2m_2}
  \nn \\
  & & \quad \stackrel{(\ref{WigPro})}{=}
  \sum_{M'} \mu(e^{J}_{M'}(j,j)) \otimes T^{(J)}_{M'M}.
  \nn
\eea

 Employing (\ref{covarmu}) we now extract a set of covariant relations 
under $\varphi_R.$ 
The $J=0$ relation 
$ \varphi_R \circ \mu(e^0_0(j,j)) = \mu(e^0_0(j,j)) $ 
signifies that $ \mu(e^0_0(j,j))$ is a scalar under the right coaction. 
It may be equated to a constant parameter $r$
\beq
  \mu(e^0_0(j,j)) 
  = \sum_{m_1,m_2} C^{\ jj\,0}_{m_1m_2\,0} e^j_{m_1} e^j_{m_2} = r.
  \label{J0cov}
\eeq
If $ J = j, $ then $ \mu(e^j_m(j,j)) $ and $ e^j_m $ transform identically 
under $\varphi_R.$ Therefore $ \mu(e^j_m(j,j))$ is, in general, 
proportional to $e^j_m$. 
It may be noted that the following relations are covariant
\beq
  \mu(e^j_m(j,j)) = \sum_{m_1,m_2} C^{\; jjj}_{m_1m_2m} e^j_{m_1} e^j_{m_2}
  = \xi e^j_m,
  \label{Jjcov}
\eeq
where the proportionality constant $ \xi \rightarrow 0 $ as $ q \rightarrow 1. $ 
For $ J \neq 0, j,$ the element $\mu(e^J_M(j,j))$ can not be a scalar, 
nor proportional to $ e^J_M$ as they transform differently. 
The relevant covariant relations are, therefore, of the form 
\beq
   \mu(e^J_M(j,j)) = \sum_{m_1,m_2} C^{\;jjJ}_{m_1m_2M} e^j_{m_1} e^j_{m_2} = 0.
   \label{Jcov}
\eeq

  As will be seen from the examples given in the next subsection, 
the simultaneous use of 
all relations from (\ref{J0cov}) to (\ref{Jcov}) gives an 
inconsistent result, since some of them do not have correct classical limits. 
In order to obtain 
a consistent covariant algebra, we have to make a choice regarding the
relations to be used for defining the algebra. Then the consistency has 
to be verified. As it is clear from the above discussion, the covariant 
algebras can have at most two more parameters ($r, \xi$) in addition to 
the deformation parameter $q.$ It is emphasised that the origin of the 
parameters is clearly explained in the framework of the representation theory. 
We have formulated a method to construct $\SL$-covariant algebras with 
respect to the right coaction. It is possible to repeat the same discussion 
for the left coaction.

%
\subsection{Quantum plane and quantum spheres}
\label{QplaneQsphere}

  We apply the general prescription in the previous subsection 
to $ j=1/2 $ and $ j=1 $ corepresentations. As will be seen, 
the obtained covariant algebras correspond to the quantum 
plane of Manin for $j=1/2$ 
and the quantum spheres of Podle\'s for $j=1.$ 

  Let us start with $ j=1/2 $ case where the relevant tensor 
product decomposition is given by 
$ 1/2 \otimes 1/2 = 1 \oplus 0. $ 
We denote the basis of $ V^{(1/2)} $ by 
$ (x,y) = (e^{1/2}_{1/2}, e^{1/2}_{-1/2}). $ 
The quantum matrix which coacts on this basis is 
given by (\ref{thalfmatrix}). 
Using explicit formula of CGC given in \cite{KS,BL}, 
we obtain from (\ref{J0cov}) for $ J = 0$
\beq
  x y - q y x = r. 
  \label{plane-r}
\eeq
If we set $ r=0,$ then (\ref{plane-r}) is reduced to the 
quantum plane relation. 
For $J=1,$ we obtain, from (\ref{Jcov}), unacceptable relations 
such as $ x^2 = y^2 = 0.$ Thus we take only (\ref{plane-r}) as 
defining relations of our covariant algebra. 

  We next investigate $ j = 1 $ case, namely, the adjoint 
corepresentation of $ SL_q(2).$ Since the adjoint corepresentation 
of the classical $ SL(2) $ corresponds to the fundamental corepresentation 
of $ SO(3),$ the covariant algebra may be interpreted as a sphere. 
The relevant tensor product decomposition is 
$ 1 \otimes 1 = 2 \oplus 1 \oplus 0.$ 
The basis of $ V^{(1)}, $ on which the quantum matrix 
(\ref{t1matrix}) coacts, is denoted by $ x_m = e^1_m. $ 
The covariant relation for $J = 0 $ is obtained via (\ref{J0cov})
\beq
  x_0^2 - qx_1 x_{-1} - q^{-1} x_{-1}x_1 = r.
  \label{q-radius}
\eeq
Explicit constructions for the $ J = 1 $ case are obtained via 
(\ref{Jjcov})
\bea
  & & (1-q^2) x_0^2 + q x_{-1} x_1 - q x_1 x_{-1} = \xi x_0,
  \nn \\
  & & x_{-1} x_0 - q^2 x_0 x_{-1} = \xi x_{-1}, 
  \quad
  x_0 x_1 - q^2 x_1 x_0 = \xi x_1.
  \label{qspheredef2}
\eea
For $ J = 2, $ we obtain, from (\ref{Jcov}), unacceptable relations 
such as $ x_{\pm 1}^2 = 0.$ 
Thus we take (\ref{q-radius}) and (\ref{qspheredef2}) as defining relations of 
our covariant algebra. 
We need to check the following conditions in order to verify 
whether or not the algebra is well-defined: 
\begin{enumerate}
\renewcommand{\labelenumi}{(\alph{enumi})}
 \item The constant $r$ commutes with all generators
 \item Product of three generators, say $ x_1 x_0 x_{-1}$, has two ways 
of reversing its ordering:
\[
  \begin{array}{ccccccc}
      &  & x_1 x_0 x_{-1} & \longrightarrow & x_1 x_{-1} x_0 & & \\
      & \nearrow & & & & \searrow &  \\
      x_0 x_1 x_{-1} & & & & & & x_{-1} x_1 x_0. \\
      & \searrow & & & & \nearrow & \\
      & & x_0 x_{-1} x_1 & \longrightarrow & x_{-1} x_0 x_1 & &
  \end{array}
\]
These two ways give the same result. 
\end{enumerate}
It is straightforward to verify that the conditions (a) and (b) are satisfied. 
The covariant algebra defined by (\ref{q-radius}) and (\ref{qspheredef2}) 
was first introduced by Podle\'s and interpreted as a noncommutative 
extension of 2-sphere\cite{Po87}. 
The constant $r$ in (\ref{q-radius}) may be regarded 
as square of radius. While the parameter $ \xi$ in (\ref{qspheredef2}), 
which goes to zero in 
the classical limit, does not exist in a commutative sphere. 
We thus obtained a one-parameter family of noncommutative 2-spheres.  

  By taking the higher dimensional corepresentations, 
one may systematically obtain 
higher dimensional noncommutative spaces covariant under $ SL_q(2).$ 

%
\section{Covariant superspaces of \boldmath{$ OSp_q(1|2) $}}
\label{CovSpeOSP}

\subsection{General prescription}

  In this section, the discussions in the preceding 
sections are extended to a quantum supergroup in order 
to obtain noncommutative superspaces. 
Since the representation theories of quantum algebra $ U_q[sl(2)] $ 
and quantum superalgebra $ U_q[osp(1|2)] $ are quite parallel, 
one can establish a prescription for constructing $OSp_q(1|2)$-covariant 
algebras similar to the one for $SL_q(2)$ by repeating the same 
discussion as \S \ref{SL2cov}. 
We give our results without proofs, 
since the results in this section have already published in\cite{AC05}. 
Readers may refer to \cite{AC05} for details. 

The universal enveloping algebra $\U = U_q[osp(1|2)]$ is generated 
by the two even $ K^{\pm 1}$, and the two odd elements $ v_{\pm} $ 
satisfying the commutation properties \cite{KR}
\bea
  & & K K^{-1} = K^{-1} K = 1, \qquad
      K v_{\pm} = q^{\pm 1/2} v_{\pm} K, \nn \\
  & & \{ v_+, v_- \} = -\frac{K^2-K^{-2}}{q^4-q^{-4}}. \label{Uqdef}
\eea
Each irreducible representation (finite dimensional) of the algebra $\U$ 
is specified by a nonnegative integer $\ell$ and the corresponding 
$ (2\ell+1) $ dimensional representation space $ V^{(\ell)} $ 
is also $ {\mathbb Z}_2$ graded. 
Let $ \{\; e^{\ell}_m(\lambda) \ | \ m = \ell, \ell-1, 
\cdots, -\ell \;\} $ be a basis of $ V^{(\ell)},$ where each basis 
vector has a definite parity. The index $ \lambda = 0, 1 $ 
specifies the parity of the highest weight vector $ e^{\ell}_{\ell}(\lambda). $ 
The parity of $ e^{\ell}_m(\lambda) $ equals $ \ell - m + \lambda,$ 
as it is obtained by the application of $ v_-^{\ell-m} $ on 
$ e^{\ell}_{\ell}(\lambda).$

  Tensor product of two irreducible representations of $\U$ has 
been discussed in \cite{KR,MM}. 
It is, in general, reducible and 
decomposed into a direct sum of irreducible representations.  
The rule of decomposition is identical to the classical case:
\beq
  \ell_1 \otimes \ell_2 = 
  \ell_1 + \ell_2 \oplus \ell_1 + \ell_2 \oplus \ell_1 + \ell_2-1 
  \oplus \cdots \oplus |\ell_1 - \ell_2|.
  \label{IrrepdecompOSP}
\eeq
The irreducible basis of the tensor product representations is obtained 
by using the CGC:
\beq
    \e{\ell}{m}(\ell_1,\ell_2,\Lambda) = \sum_{m_1,m_2} 
    \CGC{\ell_1\ \ell_2\ \;\ell}{m_1\,m_2\,m}\,\,
    \e{\ell_1}{m_1}(\lambda) \otimes \e{\ell_2}{m_2}(\lambda),
    \label{CGCdef}
\eeq
where $m= m_1+m_2, $ and 
$ \Lambda = \ell_1+\ell_2+\ell \ ({\rm mod}\ 2) $ is the parity of the highest 
weight vector $ \e{\ell}{\ell}(\ell_1,\ell_2,\Lambda). $ 
The CGC for the algebra $\U$ has been computed in \cite{AC05,MM} and 
the orthogonality relations similar to (\ref{CGCortho1}), (\ref{CGCortho2}) 
have also been obtained.

  On the contrary to $SL_q(2)$,  explicit expressions of corepresentation 
for $ \A = OSp_q(1/2)$ have not known yet. 
Employing the duality of the algebras $ \U $ and $ \A, $  
one can obtain the hitherto unknown corepresentations of $ \A $ from the 
already known irreducible representations of $ \U. $  
Let $ D^{\ell}(X;\lambda) $ be a representation matrix of $ X \in \U $ 
on $ V^{(\ell)}$ 
\beq
   X \e{\ell}{m}(\lambda) = \sum_{m'} \e{\ell}{m'}(\lambda) \, 
   \D{\ell}{m'm}(X;\lambda).
   \label{Xrep}
\eeq
We define a corepresentation matrix  $ T^{(\ell)}(\lambda) $ of $ \A $ 
via the duality relation
\beq
   \D{\ell}{m'm}(X;\lambda) = 
   (-1)^{\h{X}(\ell-m'+\lambda)}\pair{X}{\T{\ell}{m'm}(\lambda)},
   \label{DtoT}
\eeq
and the parity may be assigned as 
\beq
 \widehat{\T{\ell}{m'm}}(\lambda) = m'+m \ ({\rm mod}\ 2). 
 \label{parityT} 
\eeq
With this corepresentation matrix, one can show that $ V^{(\ell)} $ 
is a right comodule of $\A$ and that $ T^{(\ell)}(\lambda) $ 
satisfies the product law similar to (\ref{WigPro}). 
It is not difficult to find $ T^{(\ell)}(\lambda) $ for lower values of 
$ \ell $ from (\ref{DtoT}). 
For $ \ell = 1, $ we obtain
\beq
  T^{(1)}(0) = \left(
   \begin{array}{rrr}
      a & \alpha & b \\
      \gamma & e & \beta \\
      c & \delta & d
   \end{array}
  \right),
  \quad
  T^{(1)}(1) = \left(
   \begin{array}{ccc}
      a & -\alpha & b \\
      -\gamma & e & -\beta \\
      c & -\delta & d
   \end{array}
  \right),
  \label{T1}
\eeq
where the entries in latin (greek) characters are of even (odd) parity. 
For $ \ell=2,$ the entries of corepresentation matrix are quadratic in 
$ \ell = 1 $ entries
\bea
  & & T^{(2)}(0) = 
  \nn \\
  & & {\small
  \left(
     \begin{array}{ccccc}
        a^2 & \kappa_1 a \alpha & \kappa_3 ab & \kappa_1 \alpha b & b^2
        \\
        \kappa_1 a \gamma & ae + q^{-1}\gamma \alpha & \kappa_2(a\beta + q^{-1} \gamma b)
        & -\alpha \beta + q^{-1} e b & \kappa_1 b \beta
        \\
        \kappa_3 ac & \kappa_2(a \delta + c \alpha) 
        & ad + q^{-1}[2] \alpha \delta + q^{-2} b c & \kappa_2(\alpha d + \delta b) 
        & \kappa_3 bd
        \\
        \kappa_1 \gamma c & \gamma \delta + q^{-1} c e 
        & \kappa_2(\gamma d + q^{-1} c \beta) & ed + q^{-1}\beta \delta
        & \kappa_1 \beta d
        \\
        c^2 & \kappa_1 c \delta & \kappa_3 c d & \kappa_1 \delta d & d^2
     \end{array}
  \right),}
  \nn \\ 
  \label{T2mat0}
\eea
where
\beq
   \kappa_1 = \sqrt{\frac{[4]}{q[2]}}, 
   \quad
   \kappa_2 = \sqrt{q^{-1}[3]}, 
   \quad
   \kappa_3 = \kappa_1 \kappa_2,
   \quad
   [n] = \frac{q^{-n/2}-(-1)^n q^{n/2}}{q^{-1/2}+q^{1/2}}.
   \label{kappas}
\eeq

  We are ready to discuss covariant algebras for the quantum 
supergroup $\A.$  Following the arguments for $ SL_q(2),$ 
we define the composite object
\beq
  E^{L}_{M} \equiv \mu( \e{L}{M}(\ell,\ell,\Lambda) )
  = \sum_{m_1,m_2} \CGC{\;\ell\ \;\,\ell\ \;\;L}{m_1\,m_2\,M} 
  \e{\ell}{m_1}(\lambda) \e{\ell}{m_2}(\lambda),
  \label{ELM}
\eeq
where $ \Lambda = L \ ({\rm mod}\ 2). $ 
Then the right coaction on $ E^{L}_{M}$ is shown to be
\beq
   \varphi_R(E^{L}_{M}) = \sum_{M'} E^{L}_{M'} \otimes \T{L}{M'M}(\Lambda).
   \label{phiELM}
\eeq
A following set of covariant relations are extracted from (\ref{phiELM}) 
depending on the values of $L$
\bea
  & & 
  E^{0}_{0}(0) 
  = \sum_{m_1,m_2} \CGC{\;\ell\ \;\,\ell\ \;\;0}{m_1\,m_2\,0}\,\, 
  \e{\ell}{m_1}(\lambda) \e{\ell}{m_2}(\lambda)
  = r, \quad (L = 0)
  \label{E00} \\
  & &
    E^{\ell}_m(\lambda)
  = \sum_{m_1,m_2} \CGC{\;\ell\ \;\,\ell\ \;\;\ell}{m_1\,m_2\,m}\,\, 
  \e{\ell}{m_1}(\lambda) \e{\ell}{m_2}(\lambda) 
  = \xi \e{\ell}{m}(\lambda),
  \quad (L = \ell)
  \label{E1m} \\
  & & 
    E^{L}_{M} = \sum_{m_1,m_2} \CGC{\;\ell\ \;\,\ell\ \;\;L}{m_1\,m_2\,M}\,\, 
  \e{\ell}{m_1}(\lambda) \e{\ell}{m_2}(\lambda)
  = 0, 
  \quad (L \neq 0, \ell)
  \label{ELMcov}
\eea
In (\ref{E1m}), the proportionality constant $ \xi $ is of even parity 
if $ \lambda = \ell $ (mod 2), or odd parity if $ \lambda \neq \ell $ 
(mod 2). 
As already seen in the case of $ SL_q(2),$ the simultaneous use of 
all relations from (\ref{E00}) to (\ref{ELMcov}) gives an 
inconsistent result. We have to make a choice of appropriate 
relations defining a covariant algebra. Then we should check 
the consistency conditions (a) and (b) given in \S \ref{QplaneQsphere}. 

%
%
\subsection{Quantum superspace and quantum superspheres}

  As an application of the prescription given in the previous 
subsection, let us examine the covariant algebras corresponding 
to $ \ell = 1, 2 $ with $ \lambda = 0. $ 
The covariant algebra for $ \ell = 1 $ is identified with 
the quantum superspace. The one for 
$ \ell = 2 $ is interpreted as a noncommutative extension of 
supersphere.

  We start with the case of $ \ell = 1, $ where the 
relevant tensor product decomposition is given by 
$ 1 \otimes 1 = 2 \oplus 1 \oplus 0.$ 
We denote the basis of $ V^{(1)}$ by $ z_m = \e{1}{m}(0)$ on 
which the quantum supermatrix $ T^{(1)}(0) $ in (\ref{T1}) coacts.  
Thus $ z_{\pm 1} $ are 
parity even and $ z_0 $ is parity odd. 
Using the CGC given in \cite{AC05}, we obtain from (\ref{E00}) for 
$ L = 0 $
\beq
   q^{1/2}z_{-1} z_1 + z_0^2 - q^{-1/2} z_1 z_{-1} = r.
   \label{rad10}
\eeq
For $ L = 1, $ we have $ \Lambda \neq \lambda $, and, therefore, the 
parameter $\xi$ is a Grassmann number: 
\bea
   & & -q^{1/2} z_0 z_1 + q^{-1/2} z_1 z_0 = \xi z_1,
   \nn \\
   & & z_{-1} z_1 + (q^{-1/2} + q^{1/2}) z_0^2 - z_1 z_{-1} = \xi z_0,
   \label{plane1} \\
   & & q^{1/2} z_{-1} z_0 - q^{-1/2} z_0 z_{-1} = \xi z_{-1}.
   \nn
\eea
For $ L = 2, $ we obtain, using (\ref{ELMcov}), 
unacceptable relations such as $ z_1^2 = 0. $ 
We thus take (\ref{rad10}) and (\ref{plane1}) as defining relations of 
our covariant algebra. 
It is not difficult to see that the consistency condition (a) is satisfied, 
while the condition (b) requires setting $ \xi = 0. $ 
Therefore, we define our covariant algebra by combining relations 
(\ref{rad10}) and (\ref{plane1}), while maintaining $ \xi = 0$:
\bea
   & & z_1 z_0 = q z_0 z_1, \qquad z_0 z_{-1} = q z_{-1} z_0,
   \nn \\
   & & z_1 z_{-1} = q^2 z_{-1} z_1 - q(q^{-1/2}+q^{1/2}) r,
   \label{plane2} \\
   & & z_0^2 = -q^{-1}[2] z_1 z_{-1} - q^{-1}r.
   \nn 
\eea
This may be interpreted as the most general form of a quantum superspace. 
The simplest quantum superspace corresponds to the choice of $ r = 0. $

  We next investigate a covariant algebra for $ \ell = 2. $ 
This may be interpreted as a supersymmetric extension of a noncommutative 
sphere, since $ \ell = 2 $ corresponds to the adjoint representation of 
the algebra $\A.$ The quantum supersphere may have applications in 
integrable quantum field theories. Some models of integrable field 
theory which has $ osp(1|2) $ symmetry and in which supersphere appears as a 
target space have been considered\cite{SWK}. 
If an extension of such models having quantum algebra symmetry 
is considered, the quantum supersphere will 
also appear as a target space. 

 Let us denote the basis of $ V^{(2)} $ by $ Y_m = \e{2}{m}(0)$, where   
$m = 0, \pm 1, \pm 2.$ Here $ Y_{0}, Y_{\pm 2} $ are of even parity, 
and $ Y_{\pm 1} $ are of odd. 
We seek a covariant algebra under the right coaction 
of the quantum supermatrix $ T^{(2)}(0)$ in (\ref{T2mat0}). 
In order to regard the obtained covariant algebra as a noncommutative 
extension of a supersphere, we need one relation defining the radius of 
the supersphere 
and ten commutation relations of supersphere coordinates $ Y_m. $ 
In addition to those relations, two more relations which 
relate $ Y_{\pm 1}^2 $ to other coordinate are needed, 
since the odd elements may loose their nilpotency 
at the quantum level. We thus have to find thirteen relations to 
define the quantum supersphere.  

  The relation for radius is obtained via (\ref{E00}), {\it i.e.} 
$ L = 0$ 
\beq
  q^{-1} Y_2 Y_{-2} - q^{-1/2} Y_1Y_{-1} - Y_0^2 + q^{1/2} Y_{-1} Y_1 
  + q Y_{-2} Y_2 = r, \label{radius1}
\eeq
where $r$ is a constant corresponding to the square of radius. 
As commutation relations of the coordinates $ Y_m, $ 
we admit the sets of covariant relations for $ L = 2 $ and $ L = 3$ 
obtained via (\ref{E1m}) and (\ref{ELMcov}).  
Each of them  contains five and seven relations, respectively. 
We now have obtained the required number of commutation relations 
and it is easy to verify that their 
classical limit coincide with the commutative supersphere. 
To test whether they consistently define an algebra, 
we need to check for the conditions (a) and (b) mentioned in 
(\ref{QplaneQsphere}). 
It may be proved by direct computation that the said conditions 
are, however, not satisfied. 
In order to make the algebra well-defined, we incorporate the 
$ L = 1 $ relations. With the aid of $ L = 1 $ relations, 
one can verify that the consistency conditions are satisfied. 
The remaining $ L = 4 $ relations can not 
be incorporated, since they contain unacceptable 
equations such as $ Y_{\pm 2}^2 = 0. $ 

  As a result, we have sixteen relations. 
As all the relations are covariant by construction, 
their linear combinations are also covariant. 
Taking linear combinations, the relations which defined 
the quantum supersphere covariant under the coaction of 
the algebra $\A$ are summarized as follows: 
the radius relation (\ref{radius1}), ten commutation 
relations, two relations for $ Y_{\pm 1}^2 $ and 
three constraints. The classical limit of three 
constraints are not required in the commutative case. 
However, we need the constraints to make our algebra 
well-defined. 
The relations for $ Y_{\pm 1}^2 $ show that the 
odd coordinate of commutative spheres are no longer 
nilpotent in the noncommutative setting. Some of the 
defining relations contains one additional parameter 
$ \xi $ originated in (\ref{E1m}). We thus have obtained 
a one-parameter family of noncommutative superspheres. 
Explicit expressions of the defining relations of the quantum 
supersphere are found in \cite{AC05}. 

  Before closing this section, we briefly mention some 
properties of the quantum superspheres obtained above. 
They enjoy two realizations: 
The first one is the realization by 
$\U$-covariant oscillator introduced in \cite{TW}. 
This realization allows us, via realizing the covariant 
oscillator in terms of conventional $q$-oscillator, to 
obtain a infinite dimensional matrix representation of our 
quantum supersphere. 
In the second realization, the coordinates $ Y_{m} $ 
of quantum superespheres are expressed in terms of the elements of 
$\A.$ More precisely, $ Y_k$ is a linear combination of 
the entries of $k$th column of the adjoint corepresentation 
matrix (\ref{T2mat0}). Therefore, the quantum supersphere can 
be regarded as a subalgebra of $ \A.$ 
This subalgebra is specified by the infinitesimal characterization 
which was first developed for $ SL_q(2)$ \cite{DK}. 
The infinitesimal characterization tells us that 
amongst subalgebras of $ \A,$ the quantum supersphere is the 
one annihilated by a linear combination of the twisted primitive 
elements of $ \U. $ 
An element $u \in \U $ possessing a coproduct structure 
$ \Delta(u) = g \otimes u + u \otimes g^{-1}$ with $ g \in \U $ being 
a group-like element is said to be twisted primitive with respect 
to $g$. 
There exist three twisted primitive elements in $\U,$ that is,  
$ v_{\pm} $ and $ K-K^{-1}. $ 
We now define an action of an element of $ u \in \U $ on  
$ a \in \A $ by
\beq
   a \odot u = (-1)^{\h{a}\h{u}}(u \otimes {\rm id})(\Delta(a)) 
   = \sum (-1)^{\h{a}\h{u}}  \pair{u}{a_{(1)}} a_{(2)},
     \label{au}
\eeq
where Sweedler's notation for coproduct, 
$\Delta(a) = \sum a_{(1)} \otimes a_{(2)}$, is used and 
$ \h{a}, \, \h{u} $ denote the parity of the elements $ u,\, a. $ 
For a twisted primitive element $ u, $ it is straightforward to 
verify that 
\beq
   a \odot u = 0 \quad {\rm and} \quad b \odot u 
   = 0 \quad \Rightarrow \quad (ab) \odot u = 0.
   \label{aubuabu}
\eeq
Thus a set of elements of $ \A $ annihilated by a twisted primitive element $u$ form 
a subalgebra of $\A.$ 
Indeed, the quantum supersphere realized in terms of  $T^{(2)} $ in (\ref{T2mat0}) 
is a subalgebra of $ \A $ that is annihilated by the twisted primitive element 
$ {\cal P}_R$
\bea
   & & {\cal P}_R = -\sqrt{g_3}\, v_+ + \sqrt{g_1}\, v_-,
   \label{PR} \\
   & & Y_k \odot {\cal P}_R = 0, \qquad k = \pm 2,\; \pm 1,\; 0. \label{YkPR} 
\eea
$ {\cal P}_R $ consists of only odd twisted primitive elements.  
This is a difference from the quantum sphere for $SL_q(2)$. In that example, 
all the twisted primitive elements contribute to the annihilation 
operator of quantum sphere.


\section{Concluding remarks}
\label{ConcludingR}

  We have developed a common general prescription for constructing 
noncommutative covariant spaces of $ SL_q(2)$ and $ OSp_q(1|2). $ 
By this construction, it is possible to obtain covariant algebras 
for a given representation of $ SL_q(2) $ or $ OSp_q(1|2).$ 
Indeed, the known noncommutative spaces, namely Manin's quantum 
plane and Podl\'es' quantum spheres, were recovered for  $SL_q(2)$ 
and novel ones, 
{\it i.e.} their extensions to $OSp_q(1|2), $ were obtained. 
We note a difference of the present work from others\cite{Ma,Ma89,VS90}. 
In order to obtain higher dimensional noncommutative spaces, 
we use a higher dimensional representation of the fixed 
quantum group $ SL_q(2) $ (or $ OSp_q(1|2)$), 
while higher rank quantum groups are used 
in \cite{Ma,Ma89,VS90}. 

  We believe that the results of this work are useful for making 
progress in constructing supersymmetric versions of noncommutative 
geometry. For instance, we construct noncommutative superspace, say 
quantum supersphere, by our method. 
Then we may consider differential calculi on the space.  
It allows us to compute its curvature, metric and so on 
based on the framework of Ref. \cite{AC04}. 
It may also be possible to extend our method to higher rank quantum 
(super) groups by taking into account the multiplicity of irreducible 
decomposition of tensor product representations.

%


\section*{Acknowledgments}

The work of N.A. is partially supported by the 
grants-in-aid from JSPS, Japan (Contract No. 15540132). 
The other author (R.C.) is partially supported by the grant 
DAE/2001/37/12/BRNS, Government of India. 


\end{document}